\documentclass[12pt]{article} 
\usepackage{hfstyle}
\usepackage{graphicx}
\usepackage{amsmath, amssymb,bbm}
\usepackage[T1]{fontenc}
\usepackage{graphicx}

 \usepackage{mathptmx}      
\newtheorem{definition}{Definition}
\newtheorem{corollary}{Corollary}
\newtheorem{proposition}{Proposition}
 \newcommand{\NN}{{\mathbbm{N}}}

\begin{document}
\title{Adam Adamandy Kocha\'nski's approximations of $\pi$: reconstruction of the algorithm
\footnote{This is updated version (with typos corrected) of the paper which appeared in
\emph{The Mathematical Intelligencer} 34(4), 40--45 (2012).}}
\author{Henryk Fuk\'s 
      \oneaddress{
         Department of Mathematics and Statistics, Brock University\\
         St. Catharines, Ontario  L2S 3A1, Canada\\
         Toronto, Ontario M5T 3J1, Canada\\
         \email{hfuks@brocku.ca}
       }
   }

%
\Abstract{%
In his 1685 paper ``Observationes cyclometricae'' published in Acta Eruditorum, Adam Adamandy Kocha\'nski
presented an approximate ruler-and-compass construction for rectification of the circle. It is not generally
known that the first part of this paper included an interesting sequence of rational approximations
of $\pi$. Kocha\'nski gave only a partial explanation of the algorithm used to produce these approximations,
while promising to publish details at a later time, which has never happened. We reconstruct the complete
algorithm and discuss some of its properties. We also argue that Kocha\'nski was very close to
discovery of continued fractions and convergents of $\pi$.
}
\maketitle

\section{Introduction}
Adam Adamandy Kocha\'nski SJ (1631--1700) was a Polish Jesuit mathematician, 
inventor, and polymath. His interest were very diverse, including problems of geometry, 
 mechanics, and astronomy, design and construction of mechanical clocks, \emph{perpetuum mobile} and mechanical
computers, as well as many other topics. He published relatively little, and most of his mathematical
 works appeared in \emph{Acta Eruditorum} between 1682 and 1696. He left a reach correspondence, however,
which currently consists of 163 surviving letters \cite{LisiakGrzebien05}. These letters include correspondence with 
Gottfried Leibniz, Athanasius Kircher SJ, Johannes  Hevelius, Gottfried  Kirch, and many other luminaries of the 17-th century, giving a rich
record of Kocha\'nski's activities and a vivid description of the intellectual life of the period.
Recently published comprehensive monograph
 \cite{Lisiak05} gives detailed account of his life and work, and includes extensive
bibliography of the relevant literature.

Among his mathematical works, the most interesting and well-known is his paper on the rectification of the circle
and approximations of $\pi$,
published in 1685 in \emph{Acta Eruditorum} under the title \emph{Observationes cyclometricae ad facilitandam praxin
accomodatae} \cite{Kochanski1685}. Annotated English translation  of \emph{Observationes} with parallel Latin version
has been made available online by the author \cite{FuksObservationes2011}.

The paper has three distinct parts, the first one giving a sequence of
rational approximations of $\pi$. This will  be the main subject of this note, so more about this
will follow in the next section.

The second part of \emph{Observationes} is the one which is most often commented and quoted. 
There the author proposes  
an approximate solution of the problem of the  rectification of the circle, giving an elegant and simple
construction of a linear segment whose length approximates $\pi$. Figure~1 shows  this construction
\begin{figure}
 \begin{center}
     \includegraphics[scale=1.1]{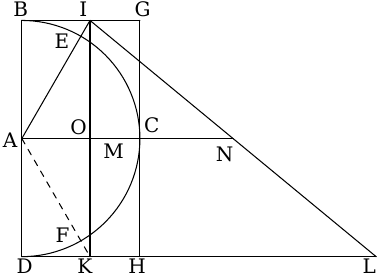}
   \end{center}
\caption{Kocha\'nski's construction of approximation of $\pi$.} 
\end{figure}
 exactly as Kocha\'nski had it the original paper. We start with drawing a semi-circle of radius
$AB=1$ centered at $A$, inscribed in a rectangle $BGHD$. Than we draw a line $AI$ such that
$\measuredangle IAC = 60^\circ$. When the lower side of the rectangle is extended
so that $HL$ is equal to the diameter of the circle, and a line is drawn from $I$ to $L$,
one can easily show that 
\begin{equation}
|IL|=\frac{1}{3}\sqrt{120-18\sqrt{3}}=3.1415333\ldots,
\end{equation}
 which agrees with $\pi$ in the first four digits after the decimal point.
This compass and ruler construction is often referred to as \emph{Kocha\'nski's construction}.

The third part of \cite{Kochanski1685} gives yet another approximation of $\pi$, this time expressing
it as a sum of multiples and fractional parts of $1/32$,
$$\frac{96}{32}+\frac{4}{32}+\frac{1}{2} \cdot \frac{1}{32} + \frac{1}{32 \cdot 32}=\frac{3217}{1024}=3.1416015625.$$

It is fair to say that if the rectification of the circle reported in \emph{Observationes} received 
a lot of attention from both contemporaries of Kocha\'nski and  
historians of mathematics \cite{Montucla1755,Cantor1880,Gunther1921},
then the first part of his paper has mostly been forgotten. In what follows we will show that 
this is perhaps unjustly so, as it
includes some intriguing sequences of fractions approximating $\pi$, origin of which has not
been explained by commentators of Kocha\'nski's work. 
\section{Sequence of rational approximations of $\pi$}
In the table on p. 395 of \cite{Kochanski1685} (also p. 2 of \cite{FuksObservationes2011}),
Kocha\'nski gives the following sequence of pairs of lower and upper rational approximants of $\pi$:
\begin{align}\label{kochanskisequence} \nonumber
\left\{\frac{25}{8},\frac{22}{7}\right\},
\left\{\frac{333}{106}, \frac{355}{113}\right\}, 
&\left\{\frac{1667438}{530762}, \frac{1667793}{530875}\right\}, 
\left\{ \frac{9252915567}{2945294501},\frac{9254583360}{2945825376} \right\}, \\ 
&\left\{\frac{136727214560643}{43521624105025}, \frac{136736469144003}{43524569930401}\right\}.
\end{align}
As mentioned in the original paper, two of these fractions can be further reduced, $\frac{1667438}{530762}=\frac{833719}{265381}$, and
$\frac{9254583360}{2945825376}=\frac{96401910}{30685681}$, while all the others are
 already written in their lowest terms. He then partially describes the  algorithm generating these approximants, 
which could be explained using modern terminology and notation as follows.
Let us denote the first element of the $n$-th pair (lower approximant) by $P_n/Q_n$, and the second element (upper approximant) by $R_n/S_n$.
The approximants are then generated by recurrence equations
\begin{align}  \label{kochconv}
  Q_{n+1}&=S_{n} x_{n}+1,\\
 P_{n+1}&=R_{n} x_{n}+3, \\
 S_{n+1}&=S_{n} (x_{n}+1)+1,\\
 R_{n+1}&=R_{n} (x_{n}+1)+3,
\end{align}
where $R_0=22$, $S_0=7$. In these formulae $x_n$ is a sequence of numbers which Kocha\'nski calls
 \emph{genitores}, giving the first four values of $x_n$:
\begin{equation}
 15,4697,5548,14774.
\end{equation}
Unfortunately, he does not explain how he obtained these numbers. He only makes the following remark regarding them:
 \begin{quote}
Methodicam pr\ae{}dictorum Numerorum Synthesin in \emph{Cogitatis,
\& Inventis Polymathematicis}, qu\ae{}, si DEUS vitam prorogaverit, utilitati
public\ae{} destinavi, plenius exponam;\footnote{I will explain the  method of generating the aforementioned numbers
 more completely in \textit{Cogitata \& Inventa 
 Polymathematica}, which work, if God prolongs my life,
  I have decided to put out for public benefit. (transl. H.F.)}
 \end{quote}
In spite of this declaration, \textit{Cogitata \& Inventa 
 Polymathematica} have never appeared in print. It is possible that some explanation
could have been found in unpublished manuscripts of Kocha\'nski, 
but unfortunately all his personal papers gathered by the National Library in Warsaw perished during the Warsaw Uprising in 1944,
when the Library was deliberately set on fire by the German occupants. 
To the knowledge of the author, nobody has ever attempted to find the algorithm for generating the sequence of
\textit{genitores}. In the most comprehensive analysis of Kocha\'nski's mathematical works published up to date \cite{Pawlikowska69}, Z. Pawlikowska  did not offer any explanation either. 

In the subsequent section, we attempt to reproduce the most likely method by which Kocha\'nski could have obtained
the sequence of \textit{genitores}, and, consequently, a sequence of approximants of $\pi$ converging to
$\pi$. We will also explain why he gave only the first four terms of the sequence. 

\section{Construction of genitores}
It seems plausible that the starting point for Kocha\'nski's considerations was Archimedes
approximation of $\pi$ by $22/7$ and the result often ascribed to
 Metius, but really due to Simon Duchesne (lat. Simonis a Quercu) \cite{Metius1626,Wielsaw2009}, 
\begin{align}
 \frac{333}{106} < \pi < \frac{355}{113}.
\end{align}
It is also likely  that Kocha\'nski then noticed
 that the fraction $\frac{333}{106}$  can be obtained from
Archimedes' approximation by using so-called \emph{proportionum intermedia}, or ``mean proportion''.
For two fractions $\frac{a}{b}$ and $\frac{c}{d}$, the mean proportion is defined as
$\frac{a+c}{b+d}$. The fraction $\frac{333}{106}$ can be written as a mean proportion of
extension of $\frac{22}{7}$ and $\frac{3}{1}$,
\begin{align}
 \frac{333}{106}&=\frac{22 \cdot 15 + 3}{7 \cdot 15 +1}.
\end{align}
Where is the factor 15  coming from? The key observation
here is that 15 is the ``optimal'' factor, in the sense that it is the largest integer value of $x$
for which    $\frac{22 \cdot x + 3}{7 \cdot x +1}$ remains smaller than $\pi$. 

The next most likely step in Kocha\'nski's reasoning was the observation
that the upper approximant can be obtained by incrementing 15 to 16,
\begin{align}
\frac{355}{113} &=\frac{22 \cdot 16 + 3}{7 \cdot 16 + 1}.
\end{align}
By repeating this procedure for $\frac{355}{113}$ one can produce another pair of
approximants,  
\begin{align}
\frac{355\cdot4697+3}{113\cdot4697+1}&=\frac{1667438}{530762}, \\
\frac{355\cdot4698+3}{113\cdot4698+1}&=\frac{1667793}{530875},
\end{align}
where $4697$ is again the largest integer $x$  for which $\frac{355 \cdot x + 3}{113 \cdot x +1}<\pi$.
Recursive application of the above process produces the desired sequence of pairs 
given in eq. (\ref{kochanskisequence}), and the values of  $x$ thus obtained 
are precisely what Kocha\'nski calls \emph{genitores}. 

What remains to be done is proving that the above algorithm indeed produces 
a sequence of lower and upper approximants of $\pi$, and that these converge
to $\pi$ in the limit of $n \to \infty$. 

\section{Kocha\'nski approximants}

We will present the problem in a general setting. In what follows,
$\alpha$ will denote a  positive irrational number which we want to approximate by rational fractions.

 Suppose that we have a pair of positive integers $R$ and $S$ such
that their ratio is close to $\alpha$ but exceeds $\alpha$,  $R/S > \alpha$. 
Together with $\lfloor\alpha\rfloor$, we then have two rational bounds on $\alpha$,
\begin{equation}\label{zerobounds}
 \frac{\lfloor\alpha\rfloor}{1} < \alpha < \frac{R}{S}.
\end{equation}

Suppose now that we want to improve these bounds. As we will shortly see, this can be achieved by considering
the mean proportion of the extension of the upper bound and the lower bound, that is,
a fraction which has the form
\begin{equation}
 \frac{Rx+\lfloor\alpha\rfloor}{Sx+1},
\end{equation}
where $x$ is some positive integer. Before we go on, let us first note 
that
the function $f(x)=\frac{Rx+\lfloor\alpha\rfloor}{Sx+1}$, treated as a function of real $x$,
has positive derivative everywhere except at $x=-1/S$, where it is undefined, and that
there  exists  $x$ where $f(x)=\alpha$, given by
$x=(\alpha-\lfloor\alpha\rfloor)/(R-\alpha S)$. 

\begin{definition} Let $\alpha$ be a positive irrational number, and let $R$ and $S$ be
positive integers such that $\frac{R}{S}> \alpha$. \emph{Genitor} of $R, S$ with respect to $\alpha$
will be defined as
\begin{equation}\label{defg}
 g_{\alpha}(R,S)=\left\lfloor \frac{\alpha-\lfloor\alpha\rfloor}{R-\alpha S} \right\rfloor.
\end{equation}
\end{definition}
Let us note that if $g_{\alpha}(R,S)$ is positive, then it  is the largest positive integer $x$ such that
$\frac{Rx+\lfloor\alpha\rfloor}{Sx+1} < \alpha$, i.e., 
\begin{equation}\label{altdef}
g_{\alpha}(R,S)=\max \{ x \in \NN : \frac{Rx+\lfloor\alpha\rfloor}{Sx+1} < \alpha\}.
\end{equation}
In this notation, the four  \emph{genitores} given in the paper can thus be written as
$g_{\pi}(22,7)=15$, $g_{\pi}(355, 113)=4697$, 
$g_{\pi}(1667793,530875)=5548$, and $g_{\pi}(9254583360,2945825376)=14774$.

Using the concept of \emph{genitores}, we can now tighten the bounds given in eq. (\ref{zerobounds}).
\begin{proposition}
 For any $\alpha \in \mathbbm{IQ}^{+}$ and $R,S \in \mathbbm{Q}^{+}$, if $\frac{R}{S}>\alpha$ and 
if the genitor $g_{\alpha}(R,S)$ is positive, then
\begin{equation}\label{propqinequalities}
\lfloor\alpha\rfloor <  \frac{R g_\alpha(R,S)+\lfloor\alpha\rfloor}{S g_\alpha(R,S)+1} < \alpha <
  \frac{R (g_\alpha(R,S)+1)+\lfloor\alpha\rfloor}{S (g_\alpha(R,S)+1)+1}<\frac{R}{S}.
\end{equation}
\end{proposition}
The second and third inequality is a simple consequence of the definition of $g_\alpha(R,S)$ 
and eq. (\ref{altdef}). The first one can be
demonstrated as follows. Since $R/S>\alpha$, then $R/S>\lfloor\alpha\rfloor$, and therefore 
 $R g_\alpha(R,S) > \lfloor\alpha\rfloor S g_\alpha(R,S)$. Now
$$ R g_\alpha(R,S)+\lfloor\alpha\rfloor > \lfloor\alpha\rfloor S g_\alpha(R,S) +\lfloor\alpha\rfloor,$$ and
$$\frac{Rg_\alpha(R,S)+\lfloor\alpha\rfloor}{Sg_\alpha(R,S)+1}>\lfloor\alpha\rfloor,$$
as required. To show the last inequality let us note that
$$%
\frac{R (g_\alpha(R,S)+1)+\lfloor\alpha\rfloor}{S (g_\alpha(R,S)+1)+1}-\frac{R}{S}
=\frac{\lfloor\alpha\rfloor S-R}{S(S (g_\alpha(R,S)+1)+1)}.
$$
Since $R/S>\lfloor\alpha\rfloor$, the numerator is negative, and the last inequality
of (\ref{propqinequalities}) follows. $\Box$

The above proposition gives us a method to tighten the bounds of (\ref{zerobounds}), and the next
logical step is to apply this proposition recursively. 

\begin{definition}\label{kochapproximants}
Let $\alpha \in \mathbbm{IQ}^{+}$ and let $R_0,S_0$ be positive integers such that $R_0/S_0>\alpha$
and $g_\alpha(R_0,S_0)>0$. \emph{Kocha\'nski approximants} of $\alpha$ starting from $R_0,S_0$ are sequences
 of rational numbers
$\{P_n/Q_n\}_{n=1}^\infty$ and $\{R_n/S_n\}_{n=0}^\infty$ defined recursively for $n\in \NN \cup \{0\}$ by
\begin{align}  \label{mainrec} 
P_{n+1}&=R_{n} x_{n}+\lfloor\alpha\rfloor, \\ \nonumber
Q_{n+1}&=S_{n} x_{n}+1,\\ \nonumber
R_{n+1}&=R_{n} (x_{n}+1)+\lfloor\alpha\rfloor,\\ \nonumber
S_{n+1}&=S_{n} (x_{n}+1)+1, \nonumber
\end{align}
where $x_n= g_\alpha(R_{n},S_{n})$. Elements of the sequence $\{P_n/Q_n\}_{n=1}^\infty$ will
be called \emph{lower approximants}, and element of the sequence $\{R_n/S_n\}_{n=0}^\infty$ 
-- \emph{upper approximants}.
\end{definition}

Note that 
\begin{align}
 P_n&=R_n-R_{n-1},\\
 Q_n&=S_n-S_{n-1},
\end{align}
therefore it is sufficient to consider sequences of $R_n$ and $S_n$ only, as these two sequences
uniquely define both upper and lower approximants.

\begin{proposition}
Kocha\'nski approximants have the following properties:
\begin{enumerate}
\item[(i)] $x_n$ is non-decreasing sequence of positive numbers,
\item[(ii)] $\displaystyle \lfloor\alpha\rfloor < \frac{P_n}{Q_n} < \alpha < \frac{R_n}{S_n}<\frac{R_0}{S_0}$ for all $n\geq1$,
\item[(iii)] $\displaystyle \frac{R_n}{S_n}$ is decreasing,
\item[(iv)] $\displaystyle \frac{P_n}{Q_n}$ is increasing,
\item[(v)] $\displaystyle \lim_{n \to \infty} \frac{R_n}{S_n}= \lim_{n \to \infty} 
\frac{P_n}{Q_n}  = \alpha$.
\end{enumerate}
\end{proposition}

For (i), because of the definition of $x_n=g_\alpha(R_n,S_n)$ shown in  eq. (\ref{defg}), we need to demonstrate that
$R_n-\alpha S_n$ is non-increasing. To do this, let us check the sign of 
\begin{align*}
 R_n - \alpha S_n -(R_{n+1}-\alpha S_{n+1})&= R_n - \alpha S_n - (R_{n}(x_n+1) + \lfloor\alpha\rfloor)\\ + 
\alpha (S_{n}(x_n+1)+1)
&=\alpha - \lfloor\alpha\rfloor - (R_n-\alpha S_n)x_n\\
&=\alpha - \lfloor\alpha\rfloor - 
(R_n-\alpha S_n)\left\lfloor \frac{\alpha-\lfloor\alpha\rfloor}{R_n-\alpha S_n} \right\rfloor.
\end{align*}
The last expression, by the definition of the floor operator, must be non-negative,
thus $R_n-\alpha S_n$ is non-increasing, and $x_n$ is non-decreasing as a result.
Now, since the definition of
Kocha\'nski approximants requires that $x_0$ is positive, all other $x_n$ must be positive too.

Property (ii) is just a consequence of the Proposition 1, which becomes clear once we note that
$$\frac{R_n}{S_n}=\frac{R_{n-1}(x_{n-1}+1) + 3}{S_{n-1}(x_{n-1}+1)+1},$$
and 
$$\frac{P_n}{Q_n}=\frac{R_{n-1}x_{n-1} + 3}{S_{n-1}x_{n-1}+1},$$
where $x_{n-1}=g_\alpha(R_{n-1}, S_{n-1})$.

To show (iii), let us compute the difference between two consecutive terms of the sequence $R_n/S_n$,
$$\frac{R_n}{S_n}-\frac{R_{n-1}}{S_{n-1}}=\frac{R_{n-1}y_{n-1}+\lfloor\alpha\rfloor}{S_{n-1}y_{n-1}+1}-\frac{R_{n-1}}{S_{n-1}},$$
where we defined $y_{n-1}=g_\alpha(R_{n-1}, S_{n-1})+1$. This yields
\begin{align*}
\frac{R_n}{S_n}-\frac{R_{n-1}}{S_{n-1}}&= \frac{(R_{n-1} y_{n-1} +\lfloor\alpha\rfloor) S_{n-1} 
- R_{n-1} (S_{n-1} y_{n-1} +1)}{S_{n-1} (S_{n-1} y_{n-1} + 1)}\\
&=\frac{\lfloor\alpha\rfloor S_{n-1}-R_{n-1}}{S_{n-1} (S_{n-1} y_{n-1} + 1)}< 0,
\end{align*}
because, by (ii), $R_{n-1}/S_{n-1}>\lfloor\alpha\rfloor$. The sequence $R_n/S_n$ is thus decreasing.
Proof of (iv) is similar and will not be presented here. 

Le us now note  that $R_n/S_n$ is bounded from below by $\alpha$ and
decreasing, thus it must have a limit. Similarly, $\frac{P_n}{Q_n}=\frac{R_n-R_{n-1}}{S_n-S_{n-1}}$ is bounded from above
by $\alpha$ and increasing, so again it must have a limit. To demonstrate (v), it is therefore sufficient to show 
that limits of $\frac{R_n}{S_n}$ and $\frac{R_n-R_{n-1}}{S_n-S_{n-1}}$ are the same, or, what is
equivalent, that 
\begin{equation}
  \lim_{n \to \infty} \left( \frac{R_n}{S_n}- \frac{R_n-R_{n-1}}{S_n-S_{n-1}}\right)=0.
\end{equation}
We start by defining $\gamma_n=R_n/S_n - P_n/Q_n$ and observing that
$$
\gamma_n=\frac{R_n}{S_n}- \frac{R_n-R_{n-1}}{S_n-S_{n-1}}=\frac{R_{n-1}S_n-R_n S_{n-1}}{S_n (S_n - S_{n-1})}.
$$
By substituting  
$R_n=R_{n-1}(x_{n-1}+1) + \lfloor\alpha\rfloor$ and
 $S_n=S_{n-1}(x_{n-1}+1)+1$, one obtains after simplification
\begin{equation}
\gamma_n=\frac{R_{n-1}-\lfloor\alpha\rfloor S_{n-1}}{S_n (S_n - S_{n-1})}=\frac{ \frac{R_{n-1}}{S_{n-1}}-\lfloor\alpha\rfloor}
 {S_n (\frac{S_n}{S_{n-1}} - 1)}.
\end{equation}
Since $R_n/S_n$ is decreasing, and starts from $R_0/S_0$, we can write 
\begin{equation}
 \gamma_n< \frac{ \frac{R_0}{S_0}-\lfloor\alpha\rfloor} {S_n (\frac{S_n}{S_{n-1}} - 1)}
=\frac{\frac{R_0}{S_0}-\lfloor\alpha\rfloor}{S_n(x_{n-1}+\frac{1}{S_{n-1}})},
\end{equation}
where we used the fact that $S_{n}=S_{n-1}(x_{n-1}+1)+1$ and where $x_{n-1}=g_\alpha(R_{n-1}, S_{n-1})$.
Since $x_n$ is non-decreasing, and $S_n$ increases with $n$, we conclude that $\gamma_n \to 0 $ ad $n\to \infty$, 
as required. $\Box$

Let us remark here that $x_n$ is indeed only non-decreasing, and it is possible for two consecutive 
values of $x_n$ to be the same.
For example, for $\alpha=\sqrt{2}$ and $R_0/S_0=3/2$, we obtain 
$$\{x_n\}_{n=0}^\infty=2, 4, 4, 15, 17, 77, 101, 119, \ldots,$$
where $x_1=x_2$.
 
\section{Initial values}
One last thing to explain is the choice of the starting values $R_0$, $S_0$. 
Definition \ref{kochapproximants} requires that the \emph{genitor} of these initial values is positive,
so how can we choose $R_0,S_0$ to ensure this?
We start by noticing that the second pair of Kocha\'nski's approximants
(fractions  $\frac{ 333}{106}$ and $\frac{ 355}{113}$) are known
to appear in the sequence of convergents of the continuous fraction  representation of $\pi$.
As we shall see, this is not just a coincidence.

Let us first recall two basic properties of continuous fraction expansion of a positive
irrational number $\alpha$,
$$
\alpha= a_0 + \cfrac{1}{a_1 + \cfrac{1}{a_2 + \cfrac{1}{a_3 + \cfrac{1}{a_4 + \ddots}}}}.
$$
By convergent $p_n/q_n$ we will mean a fraction (written in its lowest terms) obtained
by truncation of the above infinite continued fraction after $a_n$.  The first 
property we need is the recursive algorithm for generating convergents and values of $a_n$.
\begin{proposition}\label{defconvergents}
Consecutive convergents $p_n/q_n$ of $\alpha$ can be obtained by applying the recursive formula
\begin{align}  
  a_{n+1}&=\left\lfloor \frac{\alpha q_{n-1} - p_{n-1}}{p_n-\alpha q_n} \right\rfloor, \label{acf}\\
  p_{n+1}&=p_n a_{n+1} + p_{n-1},\\
  q_{n+1}&=q_n a_{n+1} + q_{n-1},
\end{align}
with initial conditions $a_0=\lfloor \alpha \rfloor$, $a_1=\lfloor \frac{1}{\alpha - a_0} \rfloor$,
 $p_0=a_0$, $q_0=1$, $p_1=a_0 a_1+1$, $q_1=a_1$. 
\end{proposition}
For example, for $\alpha=\pi$ we obtain 
\begin{multline} 
\left\{ \frac{p_n}{q_n} \right\}_{n=0}^\infty = 
\left\{ \frac{3}{1}, \frac{22}{7}, \frac{ 333}{106}, \frac{ 355}{113}, \frac{ 103993}{33102}, 
\frac{ 104348}{33215}, \frac{ 208341}{66317}, \frac{ 312689}{99532}, \right.\\ \left.  
\frac{ 833719}{265381}, \frac{ 1146408}{364913}, 
\frac{ 4272943}{1360120},
 \ldots \right\}
\end{multline}

Convergents are know to be the best rational approximations of irrational numbers, which
can formally be stated as follows.
\begin{proposition}
 If $p_n/q_n$ is a convergent for an irrational number $\alpha$ and $p/q$ is an arbitrary fraction
with $q<q_{n+1}$, then 
\begin{equation}\label{bestrational}
|q_n \alpha - p_n|< |q a -p|
\end{equation}
\end{proposition}
Elementary proofs of both of the above propositions can be found in \cite{Wall67}.
We also need to recall that convergents $p_n/q_n$ are alternatively above and below $\alpha$,
so that for  odd $n$ we always have $p_n/q_n>\alpha$, and for even $n$, $p_n/q_n<\alpha$.
Suppose that we now take some odd convergent $p_{2k+1}/q_{2k+1}$, and further set $p=\lfloor \alpha \rfloor$,
$q=1$. Inequality (\ref{bestrational}) then becomes
\begin{equation}
 p_{2k+1}-\alpha q_{2k+1}<\alpha - \lfloor \alpha \rfloor,
\end{equation}
and hence 
\begin{equation}
\frac{\alpha -\lfloor \alpha \rfloor}{p_{2k+1}-\alpha q_{2k+1}}>1.
\end{equation}
This, by the definition of the \emph{genitor} given in eq. (\ref{defg}), yields $g_\alpha(p_{2k+1},q_{2k+1})>0$,
leading to the following corollary.
\begin{corollary}
 If $p_{2k+1}/q_{2k+1}$ is an odd convergent of $\alpha$, then $g_\alpha(p_{2k+1},q_{2k+1})>0$,
and $R_0=p_{2k+1}$, $S_0=q_{2k+1}$ can be used as initial values in the construction of
Kocha\'nski's approximants. In particular, one can generate Kocha\'nski's approximants 
starting from the first convergent of $\alpha$, by taking $R_0=a_0 a_1+1$, $S_0=a_1$, where
$a_0=\lfloor \alpha \rfloor$, $a_1=\lfloor 1/(\alpha - a_0) \rfloor$.
\end{corollary}
Note that Kocha\'nski in his paper indeed started from the first convergent of $\pi$,
by taking $R_0/S_0=22/7$. 
Obviously, if one starts from the first convergent $R_0=p_1$, $S_0=q_1$, then 
the first lower approximant will be the second convergent,  $P_1=p_2$, $Q_1=q_2$, and
indeed in Kocha\'nski's case $P_1/Q_1=p_2/q_2=333/106$. 
Other approximants do not have to be convergents, and they normally aren't, 
although convergents may occasionally appear in the sequence of lower or upper
approximants. 
For example, in the case of $\alpha=\pi$, $R_1/S_1=p_3/q_3=355/113$
and $P_2/Q_2=p_8/q_8=833719/265381$.

We should also add here that the choice of the first convergent as the starting point is the most natural one.
Among all pairs $R_0,S_0$ where $S_0<106$, the only cases for which 
$g_\pi(R_0,S_0)>0$ are $R=22k$, $R=7k$, where $k\in \{1,2,\ldots,15\}$. If one wants
to obtain fractions expressed by as small integers as possible, then taking $k=1$ is an
obvious choice.


\section{Concluding remarks}
We have reconstructed the algorithm for construction  of rationals approximating
 $\pi$ used in \cite{Kochanski1685}, and we have demonstrated that it can be generalized to produce approximants
of arbitrary irrational number $\alpha$. Under a suitable choice of initial values, approximants 
converge to $\alpha$.  

Using  these results, we can generate more terms of the sequence of \emph{genitores} for $\alpha=\pi$, 
$R_0/S_0=22/7$, going beyond first four terms
found in Kocha\'nski's paper:
\begin{multline}
\{x_n\}_{n=0}^\infty=\{g_\pi(R_n,S_n)\}_{n=0}^\infty=
\{15, 4697, 5548, 14774, 33696, 61072, 111231,\\ 115985, 173819, 563316, 606004,\ldots\}.
\end{multline}
We propose to call this sequence \emph{Kocha\'nski sequence}. It has been submitted  
to the Online Encyclopedia of Integer  Sequences as A191642 \cite{A191642}, and its entry
in the Encyclopedia includes Maple and Pari code for generating consecutive terms.

Knowing that $x_n=\left\lfloor \frac{\alpha-\lfloor\alpha\rfloor}{R_n-\alpha S_n} \right\rfloor$, we can also
understand why only four terms of the sequence are given in the paper. In order to 
compute $x_n$, one needs to know $\pi$ with sufficient accuracy. For example,
20 digits after the decimal point are needed in order to compute  $x_0$ to $x_3$.
Kocha\'nski was familiar with the work of Ludolph van Ceulen, who computed 35 digits of $\pi$,
and this was more than enough  to compute $x_4$. Nevertheless, Kocha\'nski in his
paper performed all computations keeping track of ``only'' 25 digits, and this
was falling just one digit short of the precision needed to compute $x_4$.

It is also interesting to notice that the recurrence equations in
Definition \ref{kochapproximants} strongly resemble recurrence equations
for convergents $p_n/q_n$ in Proposition \ref{defconvergents}. Kocha\'nski was
always adding $3$ and $1$ to the numerator and denominator in his approximants,
because, as remarked earlier, he noticed that
\begin{align}
 \frac{22 \cdot 15 + 3}{7 \cdot 15 +1}=\frac{333}{106}, \,\,\,\,\,\,\,
  \frac{22 \cdot 16 + 3}{7 \cdot 16 + 1}=\frac{355}{113}.
\end{align}
He apparently failed to notice that
\begin{equation}
  \frac{333 \cdot 1 + 22}{106 \cdot 1 + 7}=\frac{355}{113}, 
\end{equation}
that is, instead of finding the largest $x$ for which $(22x + 3)/(7x + 1)<\pi$,
one can take the last two approximants, $22/7$ and $333/106$,
and then find the largest $x$ such that 
$(333 x + 22)/(106 x + 7)>\pi$. If he had done
this he would have discovered convergents and continued fractions. His genitores would
then be $a_n$ values in the continued fraction expansion of $\pi$. In the meanwhile,
continued fractions and convergents had to wait until 1695 when 
John Wallis laid the groundwork for their theory 
in his book \emph{Opera Mathematica} \cite{Wallis1695}.

One little puzzling detail remains, however. If we look at the Definition \ref{kochapproximants},
we notice that the sequence of lower approximants $P_n/Q_n$ starts from $n=1$, not from
$n=0$, as is the case for the upper approximants $R_n/S_n$. Indeed, $P_0, Q_0$ are not needed
to start the recursion. Nevertheless, in the table of approximants given in \cite{Kochanski1685},
in the second row there is a pair of values corresponding to $n=0$, namely $P_0/Q_0=25/8$ (in the first row of the table he also 
gives the obvious bounds $3<\pi<4$). These numbers are
 not needed in any subsequent
calculation, and Kocha\'nski does not explain where do they come from. One can only speculate 
that perhaps he 
wanted the table to appear ``symmetric'', thus he entered some arbitrary fraction
approximating $\pi$ from below as  $P_0/Q_0$.


\section*{Acknowledgements}
The author acknowledges partial financial support from the Natural Sciences 
and Engineering Research Council of Canada (NSERC) in the form of a Discovery Grant.
He also wishes to thank Prof. Danuta Makowiec for  help in acquiring relevant literature,
Rev. Prof. Bogdan Lisiak SJ for prompt and helpful replies to inquires regarding  A. A. Kocha\'nski
and his works, as well as Prof. Witold Wi\k{e}s\l{}aw for correction and clarification of 
a few historical details.

\bibliography{kochanski}

\end{document}